\documentclass{article}

\usepackage{color, xcolor, colortbl}
\usepackage{graphicx,epstopdf}
\usepackage{geometry}
\usepackage{amsmath,amssymb,amsthm}
\usepackage{algorithm}
\usepackage{algorithmic}
\usepackage{bm}
\usepackage[caption=false]{subfig}
\usepackage{appendix}
\usepackage{multirow}
\usepackage{braket}
\usepackage[english]{babel}
\usepackage{hyperref}
\usepackage[capitalize]{cleveref}
\usepackage[square,numbers,sort&compress]{natbib}
\usepackage{adjustbox}
\usepackage{xspace}
\usepackage[roman]{complexity}

\newcommand{\I}{\mathrm{i}}

\newcommand{\norm}[1]{\left\lVert#1\right\rVert}

\newcommand{\op}{\text{op}}

\newcommand{\Or}{\mathcal{O}}

\newtheorem{thm}{\protect\theoremname}
\theoremstyle{plain}
\newtheorem{lem}[thm]{\protect\lemmaname}
\theoremstyle{plain}
\newtheorem{rem}[thm]{\protect\remarkname}
\theoremstyle{plain}
\newtheorem*{lem*}{\protect\lemmaname}
\theoremstyle{plain}

\theoremstyle{plain}
\newtheorem{cor}[thm]{\protect\corollaryname}

\providecommand{\definitionname}{Definition}
\providecommand{\assumptionname}{Assumption}
\providecommand{\corollaryname}{Corollary}
\providecommand{\lemmaname}{Lemma}
\providecommand{\propositionname}{Proposition}
\providecommand{\remarkname}{Remark}
\providecommand{\theoremname}{Theorem}

\usepackage{url}
\usepackage{enumitem}
\usepackage{indentfirst}

\numberwithin{equation}{section}

\begin{document}

\title{Observable Error Bounds of the Time-splitting Scheme for Quantum-Classical Molecular Dynamics}

\author{
Di Fang\thanks{%
        Department of Mathematics, Simons Institute for the Theory of Computing, and Challenge Institute for Quantum Computation, University of California, Berkeley,  CA 94720, USA (difang@berkeley.edu)},~~
Albert Tres Vilanova\thanks{%
        Department of Mathematics, University of California, Berkeley,  CA 94720, USA (albert.tres@berkeley.edu)}.
}
\date{}

\maketitle
\begin{abstract}
Quantum-classical molecular dynamics, as a partial classical limit of the full quantum Schr\"odinger equation, is a widely used framework for quantum molecular dynamics. The underlying equations are nonlinear in nature, containing a quantum part (represents the electrons) and a classical part (stands for the nuclei). An accurate simulation of the wave function typically requires a time step comparable to the rescaled Planck constant $h$, resulting in a formidable cost when $h\ll 1$. We prove an additive observable error bound of Schwartz observables for the proposed time-splitting schemes based on semiclassical analysis, which decreases as $h$ becomes smaller. Furthermore, we establish a uniform-in-$h$ observable error bound, which allows an $\Or(1)$ time step to accurately capture the physical observable regardless of the size of $h$. Numerical results verify our estimates.
\end{abstract}

\section{Introduction}
Molecular dynamics has been a long-standing topic in computational chemistry, physics and biology. The ultimate goal is to solve
the many-body Schr\"odinger equation governed by the molecular Hamiltonian describing both electrons and nuclei, direct simulation of which is known to be a formidable challenge due to the curse of dimensionality. Mixed quantum-classical approaches are thus proposed to tackle this issue. These methods treat the nuclei and electrons separately, while mathematically certain underlying tensor product ansatz is typically assumed in deriving the couple system. Such an approximation fully exploits the multi-scale structure, that is, the scale separation between the fast-varying electrons and the slowly evolved nuclei. Some mixed approaches remain fully on the quantum level, such as time-dependent self-consistent field methods \cite{SunMiller1997,tdscf1, tdscf3, tdscf2, tdscf4} and other time-dependent mean-field approximations \cite{RungeGross1984,YabanaBertsch1996,OnidaReiningRubio2002,Ullrich2011,AnFangLin2022PT,JiaAnWangEtAl2018}, while others combine together with classical or semiclassical levels, such as the Born-Oppenheimer (or adiabatic) approximation, the mean-field type treatments \cite{MeyerMiller1979,Miller2009} and trajectory based surface hopping algorithms \cite{tully1971, tully1990, tully1998mixed, LuZhou18, FangLu2018, Lasser2007, avoided_lasser, fer_las, KammererLasser2008}. In latter approaches, the nuclear degrees of freedom are described classically by a Newtonian flow with the forces acting on the nuclei as feedback of the electronic structures calculated ``on-the-fly" (for
detailed reviews, see e.g. \cite{billing2003quantum, MarxHutter2009, marx2000ab,
tully1998mixed}).  

We focus on the quantum-classical molecular dynamics (QCMD) as proposed and justified in \cite{bornemann1996quantum}, also sometimes referred to as the Ehrenfest molecular dynamics \cite{bornemann1996quantum, drukker1999basics, schutte1999singular,
tully1998mixed} whose equations of motion are given by a Newtonian flow (for the nuclei) coupled with a Schr\"odinger equation (for the electrons):
\begin{eqnarray}
\label{SN-1}
\I h \partial _{t}\psi ^{h } &=&-\frac{h^2}{2}\Delta
_{x}\psi ^{h }+V\left( x,y\left( t\right) \right) \psi ^{h }: = H^h_t \psi^h,  \quad \psi^h(0, x) = \psi^h_0(x)\\
\label{SN-2}
\ddot{y}\left( t\right)  &=&-\nabla _{y}V_{e}^{E}\left( y\left( t\right)
\right) , \quad y(0) = y_0, \quad y'(0) = v_0 ,
\end{eqnarray}%
where$\ \psi ^{h }\left( t, x\right) \in C\left( \mathbb{R}_{t}^{+}, \mathbb{R}_{x}^{d}\right) $ is the wavefunction of the electrons, 
$h $ is the rescaled Planck constant,
 $y\left(
t\right) \in C\left( \mathbb{R}_{t}^{+}\right) $ is the position of the
nucleus, $V$ is a given potential of the whole system, and 
\begin{equation*}
V_{e}^{E}=\int_{\mathbb{R}^{d}}V\left( x,y\right) \left\vert \psi ^{h
} (t,x) \right\vert ^{2}\,dx
\end{equation*}
is often called the Ehrenfest potential, given by the quantum dynamics of
the electrons, which vividly explains the spirit of the 
``on-the-fly" simulation. 

QCMD is typically seen as an extension of the
Born–Oppenheimer approximation to the time-dependent
situation that serves as a popular tool in chemical and engineering applications (see, e.g., textbooks \cite{Ullrich2011,barth2007texts,MarxHutter2009}). The basic assumption is that the masses between particles differ significantlly, which leads to the heuristics that the heavier particle can be modeled classically while the lighter ones remain quantum. Interestingly, such heuristics can be justified on a rigorous level. Known as an approximation
of full many-body quantum dynamics, it can be derived mathematically as a partial classical limit of the full molecular Schr\"odinger equation by combining the separation of the full wave function and short wave asymptotics~\cite{bornemann1996quantum}.

The numerical difficulties associated with this type of semiclassical Schr\"odinger equation 
lie in the oscillations of order $1/h$ in both time and spatial discretizations \cite{Bao:2002fy}, which needs to be resolved for an accurate simulation of the wavefunction. However, it has been observed and formally demonstrated in
\cite{Bao:2002fy}, using the Wigner analysis, that for the linear Schr\"odinger equation with scalar potentials under suitable conditions,
a time-splitting spectral method can still capture the correct {\it physical observables} -- not the wavefunction -- when the time step is much larger than $h$, although spatially one still needs to resolve $h$.  Such interesting observations have been recently justified on a rigorous level by breakthrough results via different strategies: \cite{GolseJin21} investigates the time-splitting algorithm for the von Neumann equation by measuring the Wasserstein distance between the Husimi functions of the approximate and the exact quantum density operators, and achieves an observable error bound of $\Or(\Delta t^2 + h^{1/2})$ for the second-order splitting; and more recently, a direct observable expectation comparison is considered in \cite{LasserLubich2020}, resulting in a tighter bound of $\Or(\Delta t^2 + h^{2})$.

For nonlinear Schr\"odinger equations, nevertheless, 
one still needs to resolve $h$ temporally in general, as was numerically
demonstrated in \cite{bao2003numerical}. However, for the QCMD and other Ehrenfest models, numerical evidence presented in \cite{Jin:2017bh,FangJinSparber18} suggests the existence of a temporal meshing strategy independent of $h$ for the physical observables associated to these nonlinear Schr\"odinger type system. Though heuristic arguments based on the Wigner transform have been provided in both works, no rigorous error analysis was carried out for such nonlinear problems, which is the focus of this work.

\vspace{1em}
\noindent\textbf{Contribution:}
At the continuous level, though the semiclassical limit of the Ehrenfest dynamics has been rigorously proved in \cite{Jin:2017bh} in the weak-$\ast$ sense, the distance between the Ehrenfest dynamics and its semiclassical limit has not been understood. As the \textit{first contribution}, we establish the distance of the observable expectations between such two dynamics and show that the error is $\Or(h^2)$. This Egorov-type estimate also holds when the momentum part or the potential part in the electron dynamics is absent, which serves as a building block for the numerical analysis with the time-splitting scheme.

As the \textit{second contribution} of this paper, we prove the observable error bounds of this momentum/potential type operator splitting scheme, as proposed in \cite{Jin:2017bh, FangJinSparber18}, for the QCMD model, which extends the linear result in \cite{LasserLubich2020} to a weakly \textit{nonlinear} Schr\"odinger system. Inspired by \cite{LasserLubich2020}, the proof utilizes the observable error bounds of the QCMD model on the continuous level. Different from the observable estimates obtained from the wavefunction that exhibits $h$ scaling in the denominator, this new observable error bound, taking advantage of its semiclassical limit, removes $h$ in the denominator and exhibits an additive scaling instead. Moreover, combing two types of estimates, we present a \textit{uniform} (in $h$) observable error bound of the second-order time-splitting scheme for the QCMD, which serves as our \textit{third contribution}.

\vspace{1em}
\noindent\textbf{Organization:}
The rest of this paper is organized as follows. In \cref{sec:prelim} we introduce some preliminaries of the semiclassical descriptions for the QCMD model, including the Wigner transform and Husimi functions, and two important lemmas that are used in later proofs. We then derive the observable error bounds between the QCMD and its semiclassical limit in \cref{sec:cts_ob_err}. \cref{sec:num_analysis_err_bound} discusses the time-splitting scheme for the QCMD and analyzes the numerical errors of the observable resulting from the time-splitting strategy. Finally, we provide the numerical evidence supporting our
analytical results in \cref{sec:num_results}.

\section{Preliminaries} \label{sec:prelim}

In this section we first introduce the tools to describe the semiclassical limit of QCMD and then some preliminary lemmas used in later proofs.

\subsection{The Wigner transform and the semiclassical limit}
We revisit the Wigner transform \cite{wigner1932quantum, lions1993mesures, wigner4} and the semiclassical limit of the QCMD (or the Ehrenfest dynamics) shown in \cite{Jin:2017bh}. Let us first recall that $h $-scaled Wigner transform \cite{wigner4, lions1993mesures,
markowich1993classical} associated to
any continuously parametrized family $\psi^h \equiv \{ \psi^h \}_{0\leq h \leq 1}\in L^{2}\left( \mathbb{R}^{d}\right) $ is given by
\begin{equation*}
w^{h }[ \psi^h ] \left( x,\xi \right) =\frac{1}{\left( 2\pi
\right) ^{d}}\int_{\mathbb{R}^{d}}\psi^h\left( x-\frac{h }{2}y\right)
\overline{\psi^h}\left( x+\frac{h }{2}y\right) e^{i\xi \cdot y}\, dy.
\end{equation*}%
By Plancherel's Theorem and a change of variables one easily finds
\begin{equation*}
\left\Vert w^{h }[ \psi^h ] \right\Vert _{L^{2}\left( \mathbb{R}%
^{2d}\right) }=\frac{1}{\left( 2\pi \right) ^{d/2}h ^{d}}%
\left\Vert \psi^h\right\Vert _{L^{2}\left( \mathbb{R}^{d}\right) }^{2}.
\end{equation*}%

After Wigner transforming the Schr\"odinger equation, one finds that $$w^{h }\left(t, x, \xi \right): = w^{h }[\psi^h]\left( t, x, \xi\right)$$ satisfies the following nonlocal kinetic equation (see, e.g., \cite{lions1993mesures}):
\begin{align*}
& \partial_{t} w^{h }+ \xi \cdot \nabla _{x}w^{h }+\Theta ^{h }[
V^{h }] w^{h}=0,\quad  w^{h }(0,x, \xi)=  w^{h}[\psi^h_0] , \\
& \ddot y(t) = - \int_{\mathbb{R}} \nabla_y V(x,y(t)) w^h(t, x, \xi) \, dx \, d \xi,
\quad y(0) = y_0, \quad y'(0) = v_0,
\end{align*}
where $\Theta ^{h }[V^{h }]$ is the pseudo-differential operator
\begin{equation}\label{theta}
\left(\Theta ^{h }[ V^{h }] w^{h }\right)\left(t,  x,\xi\right) =
\frac{\I}{h (2\pi)^d } \int_{\mathbb{R}^{d}}
\left( V^{h }\left( x+\frac{h }{2}z ,y(t) \right) - V^{h }\left( x-\frac{h }{2}z , y(t)\right) \right)
\widehat{w}^{h }( t, x,z) e^{iz \cdot \xi }dz
\end{equation}
with $\widehat{w}^{h}$ denoting the Fourier transformation of $w^h$ with respect to the second variable only.
It was shown in \cite{Jin:2017bh} that one can pass the limit $h \to 0_+$ and get its semiclassical limit consisting of a Liouville equation coupled with a Newtonian flow:
\begin{eqnarray}
\label{eqn:wigner1}
\partial_{t}\mu + \xi \cdot \nabla_x \mu -\nabla_x V(x,y(t)) \cdot \nabla_\xi \mu = 0, \\
\nonumber
\ddot{y}\left( t\right)  =
-\int_{\mathbb{R}^{d}}\nabla _{y}V\left( x,y(t)\right)  \mu (t, x, \xi) \,dx\,d \xi,
\end{eqnarray}
where $ w^h \stackrel{h \to 0_+}{\longrightarrow} \mu$ and the convergence is in the weak-$\ast$ sense. 
Because the measure converges in the weak-$\ast$ topology, a direct evaluation of the distance between $w^h$ and its limit $\mu$ becomes a non-trivial task. A framework based on the Wasserstein distance has been proposed for the \textit{linear} von-Neumann equation in \cite{GolseJin21}.

\subsection{The Weyl quantization and Husimi functions}
An alternative way to view the connection between the quantum dynamics and its classical analog is through the Weyl quantization. For a Schwartz funcion $a: \mathbb{R}^{2d} \to \mathbb{R}$, $a = a(x,\xi)$, we define its Weyl quantization \cite{zworski_book} to be the operator $\op(a)$ acting on $\psi \in L^2(\mathbb{R}^d)$ by the formula
\begin{equation} \label{eqn:def_weyl}
    \op(a) \psi = \frac{1}{(2\pi h)^{d}} \int_{\mathbb{R}^{2d}} a\left( \frac{x+y}{2}, \xi\right) e^{\I \xi \cdot (x-y)/h} \psi(y) \, d \xi \, dy,
\end{equation}
so that the expectation values of an operator $A = \op(a)$ can be written in terms of the Wigner transform as
\begin{equation}\label{eqn:expec_wigner}
    \langle A \rangle_{\psi}  = \bra{\psi}A\ket{\psi}  = \int_{\mathbb{R}^{2d}} a(x, \xi) w^h[\psi](x, \xi) \, dx \,d\xi.
\end{equation}
One defines the Husimi function $\sigma^h[\psi]: \mathbb{R}^{2d} \to [0, \infty)$ as
\begin{equation}\label{eqn:def_husimi}
    \sigma^h[\psi](x, \xi) = \left(w^h[g_{(0,0)}] \ast w^h[\psi]\right)(z)= \frac{1}{(2\pi h)^{d}} \lvert \braket{g_{(x,\xi)} | \psi} \rvert^2,
\end{equation}
where $g_{(x,\xi)}(\cdot)$ is the Gaussian wave packet (or the coherent state)
\[
g_{(x,\xi)}(y) =  \frac{1}{(\pi h)^{d/4}} \exp \left({-\frac{1}{2h}|y-x|^2 + \frac{\I}{h} \xi\cdot (y-x)} \right).
\]
We remark that the Husimi function of $\psi$ with $\norm{\psi} = 1$ is in fact a non-negative probability density on the phase space, and is the modulus squared of the so-called Fourier-Bros-Iagolnitzer (FBI) transform acting on $\psi$. The expectation of $A = \op(a)$ can be described by the Husimi function, 
\begin{align*}
\langle A \rangle_{\psi}  
&= \int_{\mathbb{R}^{2d}} \left(a - \frac{h}{4}\Delta a \right)(x, \xi) \, \sigma^h[\psi](x, \xi) \, dx \,d\xi + \Or{(h^2)},
\end{align*}
as detailed in \cref{lmm:exp_husimi}.

\subsection{Important Lemmas}
We review two lemmas to be used in the proof of the paper. The first lemma is the stationary phase \cite{zworski_book}, which is used in the proof of \cref{thm:cts_egorov}.
\begin{lem}[Stationary Phase Lemma] \label{lmm:stat_phase}
Assume that $ a \in C^{\infty}_{c}(\mathbb{R}^{2d}).$ Then for each positive integer N, we have

\[
\int_{\mathbb{R}^{2d}} e^{\frac{\I}{h} \braket{x, \xi}}a(x,\xi)\,dx\,d\xi = (2\pi h)^d\left(\sum_{k=0}^{N-1}\frac{h^k}{k!}\left(\frac{\braket{\nabla_x,\nabla_\xi}}{\I}\right)^ka(0,0)+\Or(h^N)\right)
\]
as $h \rightarrow 0$, where $\braket{\cdot,\cdot} $ stands for the usual inner product.
\end{lem}

\noindent The second lemma describes the approximation error of the observable expectations in terms of the Husimi function \cite[Theorem 6.11]{LasserLubich2020}, which is used in the proof of \cref{cor:egorov_wigner_husimi}.

\begin{lem}[Expectations in terms of Husimi functions]\label{lmm:exp_husimi}

Assume the function $a: \mathbb{R}^{2d} \to \mathbb{R}$ that defines the observable $A=\op(a)$ is smooth and bounded together with all its derivatives. Then, for all $\psi \in L^2(\mathbb{R})$, we have
\[
\left|\braket{A}_{\psi}-\int_{\mathbb{R}^{2d}}\left(a(x, \xi) - \frac{h}{4}\Delta a(x, \xi) \right) \sigma^h[{\psi}] (x, \xi) \,dx\, d\xi\right|\leq C h^2\norm{\psi}^2_{L^2(\mathbb{R}^d)}, \qquad 0<t\leq T,
\]
where $C\geq 0$ depends only on the derivative bounds of $a$.
\end{lem}

\section{Error between its limit in observables} \label{sec:cts_ob_err}
We quantify the distance between the QCMD and its semiclassical limit by a direct investigating on the observable expectations. When the Hamiltonian is linear and time independent, the observable bounds (also called the Egorov's theorem) is a consequence of the semiclassical expansion and commutator estimates. On an intuitive level, one has
\[
\partial_t \left(e^{\I H t/h} A e^{-\I Ht/h} \right)= e^{\I H t/h}\frac{\I}{h} [H, A] e^{-\I Ht/h},
\]
and the commutator follows the semiclassical expansion (see, e.g., \cite[Theorem 4.12]{zworski_book} and \cite[Proposition 6.2]{LasserLubich2020})
\begin{equation}\label{eqn:comm_h_a}
[H, A] = \frac{h}{\I} \op \left( \{\mathcal{H}, a\} \right)  + \Or(h^3),
\end{equation}
for any $H$ and $A$ that are Weyl-quantizations of the phase space funtions (symbols) $\mathcal{H}(x,p)$ and $a(x,p)$, respectively,
where $[H, A] = HA - AH$ is the commutator of the operators $H$ and $A$, and $\{\mathcal{H}, a\} = \nabla_p \mathcal{H} \cdot \nabla_x a  - \nabla_x \mathcal{H} \cdot \nabla_p a$ is the Poisson bracket of the corresponding functions $\mathcal{H}$ and $a$. Note that here we follow the notation of \cite{zworski_book} that is of the opposite sign comparing to that in \cite{LasserLubich2020}. More precisely, for the linear case (see, e.g., \cite[Proposition 6.3]{LasserLubich2020} and \cite[Section 11]{zworski_book}), we can estimate in terms of the operators:
\begin{align*}
      e^{\I H t/h} A e^{-\I Ht/h}  - \op(a \circ \Xi^t) 
     = &  \int_0^t \partial_s \left( e^{\I H s/h} \op(a \circ \Xi^{t-s}) e^{-\I Hs/h} \right) \, ds
     \\ 
      = & \int_0^t  e^{\I H s/h} \left(  \frac{\I}{h}[\mathcal{H}, \op(a \circ \Xi^{t-s})] - \op(\{ \mathcal{H}, a \circ \Xi^{t-s} \}) \right)  e^{-\I Hs/h}\, ds,
\end{align*}
where we use $\Xi^t$ to denote the flow governed by the Hamiltonian $\mathcal{H}$ at time $t$. Thanks to \eqref{eqn:comm_h_a}, the difference between $e^{\I H t/h} A e^{-\I Ht/h}$ and its semiclassical limit $\op(a \circ \Xi^t)$ in the operator norm from $L^2$ to $L^2$ is bounded by $\Or(h^2)$. However, when $H$ depends on the wavefunction itself, it is not sufficient to deal with linear operators.
Our goal here is also to perform the semiclassical expansion, but due to the time-dependence in the potential and nonlinearity of the problem, we consider such expansion in terms of the wavefunction (that is, the solution), instead of the operators. One certainly does not hope such results hold for arbitrary nonlinearity. Fortunately, here the observable expectations still exhibit the semiclassical behavior thanks to the structure of the QCMD system.

In order to discuss the classical counterpart of the QCMD system, we consider the following ODEs with its flow map denoted as $\Phi^t$,
\begin{equation} \label{eqn:cts_ode_flowmap}
\left\{
\begin{split}
& \dot x = \xi \\
& \dot \xi = - \nabla_x V(x,y) \\
& \dot y = v\\
& \dot v = - \nabla_y V(x,y),
\end{split}
\right.
\end{equation}%
which can be viewed as the Lagrangian description of the semiclassical limit \eqref{eqn:wigner1}. In other words, $\Phi^t$ is governed by the total Hamiltonian $\mathcal{H}_0 := \frac{|\xi|^2}{2} + \frac{|v|^2}{2} + V(x,y)$.
We are now ready to state and prove the Egorov-type theorem for the QCMD. 

\begin{thm}[Egorov's theorem for the QCMD] \label{thm:cts_egorov} Let $a(x, \xi)$ be a Schwartz function that defines the observable $A = \op(a)$, and the potential function $V$ be smooth with its derivatives of order $\geq$ 2 all bounded. Then there exist some constant $C\geq 0$ depending on the derivative bounds of $a$ and independent of $h$ and $\psi_0$ that
\[
\lvert \braket{A}_{\psi^h(t,\cdot)} - \braket{\op(a \circ \Phi^t)}_{\psi^h_0(\cdot)}  \rvert \leq Ct  h^2.
\]
\end{thm}
\begin{proof} In this proof, we drop the superscript $h$ in the wavefunction $\psi^h$ for notational simplicity. One starts by looking at the integrands of the expectations.
\begin{align*}
& \bar \psi(t,x)\,\op(a)\,\psi(t,x) - \bar \psi_0(x) \op(a\circ\Phi^t)  \psi_0(x)
= \int_0^t \frac{d}{ds} \left( \bar \psi(s) \,\op(a\circ\Phi^{t-s})\, \psi(s,x)\right) d s
\\
=& \int_0^t   \partial_s \bar \psi(s,x) \,\op(a\circ\Phi^{t-s})\, \psi(s,x) 
+  \bar \psi(s,x) \,\op(a\circ\Phi^{t-s})\, \partial_s \psi(s,x) 
- \bar \psi(s,x) \partial_t\,\op(a\circ\Phi^{t-s})  \psi(s,x) d s
\\
 = &  \int_0^t  \frac{\I h}{2} \left( -  \Delta_x \bar \psi(s,x) \,\op(a\circ\Phi^{t-s})\, \psi(s,x) 
+  \bar \psi(s,x) \,\op(a\circ\Phi^{t-s})\,  \Delta_x \psi(s,x) \right)  \, ds 
\\
& + \int_0^t   \frac{\I}{h} \left(V(x,y(s)) \bar \psi(s,x) \,\op(a\circ\Phi^{t-s})\, \psi(s,x) 
- \bar \psi(s,x) \,\op(a\circ\Phi^{t-s})\,  V(x,y(s)) \psi(s,x) \right)  \, ds 
\\
& + \int_0^t \bar \psi(s,x) \op(\{ a\circ\Phi^{t-s}, \mathcal{H}_0\})  \psi(s,x) d s
: = \int_0^t
(\text{\uppercase\expandafter{\romannumeral1}}) + (\text{\uppercase\expandafter{\romannumeral2}}) + (\text{\uppercase\expandafter{\romannumeral3}}) \, ds,
\end{align*}
where we use \eqref{SN-1} and the fact that $\partial_t\,\op(a\circ\Phi^{t-s})  = -\op(\{ a\circ\Phi^{t-s}, \mathcal{H}_0\}) $. Denote $b = a \circ \Phi^{t-s}$, one then compute the terms one by one. Starting from the momentum terms, we have
\begin{align*}
      \Delta_x \bar \psi(s,x) \,\op(b)\, \psi(s,x) 
     & = (2\pi h)^{-d} \int_{\mathbb{R}^{2d}} \Delta_x \bar \psi(s,x) b\left(\frac{x+q}{2},p\right) e^{\I p\cdot(x-q)/ h}\psi(s, q) \,dq\, dp 
\end{align*}
and it follows from the definition of Weyl quantization and integration by parts that
\begin{align*}
\bar \psi(s,x) \op(b)\Delta_x\psi(s, x) 
&= (2\pi h)^{-d} \int_{\mathbb{R}^{2d}} \bar \psi(s,x) b\left(\frac{x+q}{2},p\right) e^{\I p\cdot(x-q)/ h}\Delta_q\psi(s, q) \,dq\, dp 
\\
&=(2\pi h)^{-d} \int_{\mathbb{R}^{2d}}\bar \psi(s,x) \Delta_q \left( b\left(\frac{x+q}{2},p\right) e^{\I p\cdot(x-q)/ h}\right)\psi(s, q) \, dq \,dp
\end{align*}
Following the chain rule by changing the differentiation variable from $q$ to $x$, one has
\begin{align*}
\bar \psi(s,x) \op(b)\Delta_x\psi(s, x) 
=& (2\pi h)^{-d} \int_{ \mathbb{R}^{2d}} \bar \psi(s,x) \Delta_x \left( b\left(\frac{x+q}{2},p\right) e^{\I p\cdot(x-q)/ h}\right) \psi(s, q) \,dq\, dp 
\\
+& (2\pi h)^{-d} \int_{ \mathbb{R}^{2d}} \frac{2 \I}{  h} \bar \psi(s,x) \nabla_q b\left(\frac{x+q}{2},p\right) \cdot p e^{\I p\cdot(x-q)/ h} \psi(s, q) \,dq\, dp,
\end{align*}
where the fact that $\nabla_q {b\left(\frac{x+q}{2},p\right)} = \nabla_x {b\left(\frac{x+q}{2},p\right)}$ is also used.
Therefore, one has
\begin{align*}
    \int (\text{\uppercase\expandafter{\romannumeral1}}) \, dx
    = &
    - \frac{\I h}{2}(2\pi h)^{-d} \int_{ \mathbb{R}^{3d}} \Delta_x \bar \psi(s,x) b\left(\frac{x+q}{2},p\right) e^{\I p\cdot(x-q)/ h}\psi(s, q) \,dq\, dp \, dx
    \\
   & + \frac{\I h}{2} (2\pi h)^{-d} \int_{ \mathbb{R}^{3d}} \bar \psi(s,x) \Delta_x \left( b\left(\frac{x+q}{2},p\right) e^{\I p\cdot(x-q)/ h}\right) \psi(s, q) \,dq\, dp \, dx 
    \\
   & + (2\pi h)^{-d} \int_{ \mathbb{R}^{3d}}  \bar \psi(s,x) \nabla_q b\left(\frac{x+q}{2},p\right) \cdot p e^{\I p\cdot(x-q)/ h} \psi(s, q) \,dq\, dp \,dx \\
   = &  \int_{\mathbb{R}^d} \bar\psi(s,x) \op(\nabla_x b \cdot p) \psi(s,x) \, dx,
\end{align*}
where the first two terms cancel thanks to the integration by parts in $x$. Note that there is no approximation error in the momentum terms and this is essentially because the Laplacian is a quadratic observable. 

Next we compute the potential terms. Since $V(x,y)$ is not necessarily quadratic, one expects some approximation errors entering the calculations.
\begin{align*}
    & (\text{\uppercase\expandafter{\romannumeral2}}) 
    = \frac{\I}{h}\left( V(x,y(s)) \bar \psi(s,x) \,\op(b)\, \psi(s,x) - \bar \psi(s,x) \,\op(b)\,  V(x,y(s)) \psi(s,x)  \right)
    \\
    = & \frac{\I}{h} (2\pi h)^{-d} \int_{ \mathbb{R}^{2d}} b\left(\frac{x+q}{2},p\right) \left(V(x, y(s)) -V(q, y(s))\right)e^{\I p\cdot(x-q)/ h} \bar \psi(s,x) \psi(s, q) \,dq \,dp \\ 
    = & - \frac{\I}{h} (2\pi h)^{-d} \int_{ \mathbb{R}^{2d}} b\left(x-\frac{q}{2},p\right) \left(V(x, y(s)) -V(x-q, y(s))\right)e^{\I p\cdot q / h} \bar \psi(s,x) \psi(s, x - q) \,dq \,dp.
\end{align*}
Apply \cref{lmm:stat_phase}, one has
\begin{align} 
       - (\text{\uppercase\expandafter{\romannumeral2}})  
    =&  \bar \psi(s,x) \langle \nabla_q, \nabla_p\rangle \left( b\left(x-\frac{q}{2},p\right) \left(V(x, y(s)) -V(x-q, y(s))\right) \psi(s, x - q) \right)|_{q = p = 0} 
    \nonumber\\
     +& \bar \psi(s,x)   \frac{h}{2\I} \langle \nabla_q, \nabla_p\rangle^2 \left( b\left(x-\frac{q}{2},p\right) \left(V(x, y(s)) -V(x-q, y(s))\right)  \psi(s, x - q) \right)|_{q = p = 0} 
    + \Or(h^2)
    \nonumber\\
    =& \bar \psi(s,x)   \nabla_p b(x,0) \cdot \nabla_x V(x, y(s))   \psi(s, x )
    + \Or(h^2) 
    \nonumber\\
     +&\bar \psi(s,x)   \frac{h}{2\I} \langle \nabla_q, \nabla_p\rangle^2 \left( b\left(x-\frac{q}{2},p\right) \left(V(x, y(s)) -V(x-q, y(s))\right) \psi(s, x - q) \right)|_{q = p = 0}.
    \label{eqn:roman1}
\end{align}
On the other hand, 
\begin{align*}
&\op(\nabla_p b\cdot \nabla_x V) \psi(s, x) \\
= &  (2\pi h)^{-d}  \int_{ \mathbb{R}^{2d}} \nabla_p b\left(\frac{x+q}{2},p\right) \,\nabla_x V\left(\frac{x+q}{2}, y(s)\right) e^{\I p\cdot(x-q)/ h}\psi(s, q) \, dq dp
\\
= &  (2\pi h)^{-d}  \int_{ \mathbb{R}^{2d}} \nabla_p b\left(x-\frac{q}{2},p\right) \,\nabla_x V\left(x-\frac{q}{2}, y(s)\right) e^{\I p\cdot q / h}\psi(s, x-q) \, dq dp 
\\
= & \nabla_p b(x,0) \,\nabla_x V(x, y(s)) \psi(s, x)
\\  + & \frac{h}{\I}  \nabla_q \cdot \left(\nabla_p \nabla_p b\left(x-\frac{q}{2},p\right) \,\nabla_x V\left(x-\frac{q}{2}, y(s)\right) \psi(s, x-q)\right)|_{q = p = 0}
+ \Or(h^2),
\end{align*}
where the stationary phase lemma is applied in the last line. One can check that $\bar \psi(s,x) \op(\nabla_p b\cdot \nabla_x V) \psi(s, x)$ and $(\text{\uppercase\expandafter{\romannumeral2}})$ share the same $\Or(1)$ and $\Or(h)$ terms, and hence 
\begin{equation} \label{eqn:roman2}
(\text{\uppercase\expandafter{\romannumeral2}})  
= - \bar \psi(s,x) \op(\nabla_p b\cdot \nabla_x V) \psi(s, x) + \Or(h^2).
\end{equation}
Combing \eqref{eqn:roman1} and \eqref{eqn:roman2}, we have
\begin{align*}
    &\int_\mathbb{R} 
(\text{\uppercase\expandafter{\romannumeral1}}) + (\text{\uppercase\expandafter{\romannumeral2}}) + (\text{\uppercase\expandafter{\romannumeral3}}) \, dx
\\
 =& \int_{\mathbb{R}^d} \bar\psi(s,x) \left( \op(\nabla_x b \cdot p) 
 - \op(\nabla_p b\cdot \nabla_x V) 
 + \op(\{b, \mathcal{H}_0\}) \right) \psi(s,x) \,dx + \Or(h^2) 
 = \Or(h^2).
\end{align*}
Integrating both sides, we get the desired result.

\end{proof}

\begin{rem}
\begin{enumerate}
    \item It follows from the proof that the a similar observable bound also holds the Schr\"odinger equation with only the kinetic part or potential part. In particular, the $\Or(h^2)$ error only comes from the potential part and there is no asymptotic error in the kinetic part.
    \item Note that the proof also works with minor modifications for linear time-dependent Hamiltonian $$-\frac{h^2}{2} \Delta + V(t,x),$$
provided that $V$ is smooth with all derivatives bounded, as well as the controlled Hamiltonian in quantum control problems \cite{AnEtAl2021, Tannor1992, ZhuRabitz1998, MadayTurinici2003}
$$-\frac{h^2}{2} u_1(t) \Delta + u_2(t) V(x),$$
where $u_1, u_2 \in [0,1]$ are the control functions and $V(x)$ smooth with all derivatives bounded. But note that for these linear cases, it can be more convenient to directly work with operators. Specifically, the Hamiltonian corresponds to the Weyl quantization of time-dependent symbols, and the Egorov's theorem in such cases can be found in, e.g., \cite[Section 11]{zworski_book}.
    \item The proof only works for Schwartz functions $a(x,\xi)$. An extension to more general symbol class $S_\delta(m)$ (see \cite[Ch. 11]{zworski_book}) may be possible and is left for future study.
\end{enumerate}
\end{rem}

\cref{thm:cts_egorov}, together with \cref{eqn:expec_wigner} and \cref{lmm:exp_husimi} provides the Egorov-type results in terms of the Wigner and Husimi functions.
\begin{cor}[Egorov's Theorem in the Wigner functions and Husimi functions]\label{cor:egorov_wigner_husimi}
Under the conditions of \cref{thm:cts_egorov}, the observable expectations in terms of the Wigner and Husimi functions satisfy
\begin{align} \label{eqn:egorov_wigner}
    \left| \langle A \rangle_{\psi^h(t, \cdot)} - \int_{\mathbb{R}^{2d}} a \circ \Phi^t(x, \xi) w^h[\psi^h_0] (x, \xi) \, dx \, d\xi \right| \leq C t h^2,
    \\  \label{eqn:egorov_husimi}
    \left|  \langle A \rangle_{\psi^h(t, \cdot)} - \int_{\mathbb{R}^{2d}} \left(a -\frac{h}{4}\Delta a\right)\circ \Phi^t(x, \xi)  \sigma^h[\psi^h_0]  (x, \xi)\, dx \, d\xi \right| \leq C t h^2,
\end{align}
respectively, where $C\geq 0$ is some constant depending on the derivative of $a$ and $V$ but independent of $h$ and $\psi^h_0$.
\end{cor}

\section{Numerical Schemes and Observable Errors}\label{sec:num_analysis_err_bound}
\subsection{Time-splitting Schemes}

The time-splitting algorithm is proposed as follows: From time $t=t_{n}=n\Delta t$ to $t=t_{n+1}=\left( n+1\right)
\Delta t$, with $\Delta t$ given, the QCMD is split into two subsystems. One contains the kinetic contributions
\begin{equation}\label{step1}
\left\{
\begin{split}
& ih \partial_{t}\psi ^{h }  =-\frac{h ^{2}}{2}\Delta _{x}\psi
^{h },   \\
& \dot y = v,\\
& \dot v = 0,
\end{split}
\right.
\end{equation}%
whose unitary propagation of the wavefunction from time $t = a$ to $t = b$ is denoted as $U_T(a,b)$, and the numerical flow in this step is denoted as $\Phi_T^{b-a}$ from time $t = a$ to time $t = b$. The other is associated with the potential contributions
\begin{equation}\label{step2}
\left\{
\begin{split}
& ih \partial_{t}\psi ^{h }= V \left( x, y(t)\right) \psi^h
,  \\
& \dot y = 0,\\
& \dot v = - \int_{\mathbb{R}^d}\nabla_y V(x,y) |\psi^h(t,x)|^2 \, dx
\end{split}
\right.
\end{equation}
and we denote the unitary evolution of the wavefunction following this potential state from time $t = a$ to $t = b$ as $U_V(a,b)$, and the numerical flow in this step is denoted as $\Phi_V^{b-a}$ from time $t = a$ to time $t = b$. Note that in the evolution $U_V(a,b)$, $|\psi^h(t,x)|^2 $ does not change with respect to time. This observation is crucially important in the consideration of its classical counterpart, which makes the resulting ODE an autonomous one.
A first order Lie splitting then corresponds to evolving \eqref{step1} from time $t_n$ to $t_{n+1}$, following by the evolution of \eqref{step2} from time $t_n$ to $t_{n+1}$. Here, we focus on the discussion of the Strang splitting given as
\begin{align} \label{eqn:strang_ehrenfest}
    & \psi^{h}_{n+1}  = U_V(t_{n}+\Delta t/2, t_{n+1}) U_T(t_n, t_{n+1}) U_V (t_n, t_{n}+\Delta t/2)\psi^{h}_{n},
    \\
    & (y^{n+1}, v^{n+1})  =  \Sigma_V^{\Delta t/2} \circ \Sigma_T^{\Delta t} \circ \Sigma_V^{\Delta t/2} (y^n, v^n) : = \Sigma_{SV}^{\Delta t}(y^n, v^n),
\end{align}
where the scheme for the Newtonian part $(y,v)$ is in fact the St\"ormer-Verlet integrator (denote as SV), which will be detailed in section \ref{sec:ob_err_strang}.
To be more precise, the second-order time-splitting from time $t_n$ to $t_{n+1}$ reads
\begin{align*} \label{eqn:strang_ehrenfest_precise}
    & \psi^{h}_{\ast} = e^{-\I \Delta t V(x, y^n)/(2h)} \psi^{h}_{n}, \quad y^\ast = y^n, \quad v^\ast = v^n - \frac{\Delta t}{2} \int \nabla_y V(x, y^n) |\psi^{h}_{n}(x)|^2 \,dx,
    \\
    & \psi^{h}_{\ast \ast} = e^{\I h \Delta t \Delta_x / 2} \psi^{h}_{\ast}, \quad  y^{\ast \ast} = y^\ast + v^\ast \Delta t , \quad v^{\ast \ast} = v^\ast,
    \\
    & \psi^{h}_{n+1} = e^{-\I \Delta t V(x, y^{\ast \ast})/(2h)} \psi^{h}_{\ast\ast}, \quad y^{n+1} = y^{\ast \ast}, \quad v^{n+1} = v^{\ast \ast} - \frac{\Delta t}{2} \int \nabla_y V(x, y^{\ast \ast}) |\psi^{h}_{\ast \ast}(x)|^2 \,dx,
\end{align*}
where one can use Fourier spectral methods to implement the exponentiation of the Laplacian and trapezoidal rules for the numerical quadrature. We focus on the time-splitting strategy (in the semi-discrete set-up) and a detailed discussion on the spatial discretization is rather standard (see, e.g., \cite{Bao:2002fy, bao2003numerical, FangJinSparber18}) and beyond the scope of this paper. Nevertheless, we point out that when employing the spectral discretization spatially,  one gets an unconditionally stable scheme because the mass of the system can be shown to be conserved, following the strategies as proposed in \cite{Bao:2002fy,FangJinSparber18}, which we do not detail here.

\subsection{Observable Errors of the Strang Splitting} \label{sec:ob_err_strang}

In this section, we aim to provide a direct error bound of the observable for the Strang splitting. The strategy is to convert the error estimate of the microscopic solver into an estimate of three parts: the errors of the macroscopic numerical solver, and the asymptotic errors of both the continuous and discrete multi-scale limit, as in commonly used in the error analysis of the asymptotic preserving schemes (see, e.g., \cite{FilbetJin10,HuJinLi17,LiuWangZhou18}) and asymptotically compatible schemes \cite{TianDu20,DEliaDuEtAl20}. Note that the corresponding classical counterparts of \eqref{step1} and \eqref{step2} 
are given as
\begin{equation}
\left\{
\begin{split}
& \dot x = \xi \\
& \dot \xi = 0  \\
& \dot y = v\\
& \dot v = 0
\end{split}
\right.
\qquad \text{and} \qquad
\left\{
\begin{split}
& \dot x = 0 \\
& \dot \xi = - \nabla_x V(x,y) \\
& \dot y = 0\\
& \dot v = - \nabla_y V(x,y),
\end{split}
\right.
\end{equation}%
respectively and denote the corresponding flow map from time $t = a$ to time $t = b$ as $\Phi_T^{b-a}$ and $\Phi_V^{b-a}$. We further define
\[
\Phi_{SV}^{\Delta t}: = \Phi_V^{\Delta t/2} \circ  \Phi_T^{\Delta t} \circ \Phi_V^{\Delta t/2}.
\]
Therefore, the Strang splitting of the QCMD as proposed in \eqref{eqn:strang_ehrenfest} corresponds to the second-order St\"ormer–Verlet integrator of the Hamiltonian system governed by the Hamiltonian $h_\text{total} = \frac{|\xi|^2}{2} + \frac{|v|^2}{2} + V(x,y)$, namely,
\begin{align*}
    & x^\ast = x^n, \quad \xi^\ast = \xi^n - \frac{\Delta t}{2} \nabla_x V(x^n, y^n), \quad y^\ast = y^n, \quad v^\ast = v^n - \frac{\Delta t}{2} \nabla_y V(x^n, y^n),
    \\
    & x^{\ast \ast} = x^\ast + \Delta t \xi^\ast, \quad \xi^{\ast \ast} = \xi^\ast,\quad  y^{\ast \ast} = y^\ast + v^\ast \Delta t , \quad v^{\ast \ast} = v^\ast,
    \\
    & x^{n+1} = x^{\ast \ast}, \quad \xi^{n+1} = \xi^{\ast \ast} - \frac{\Delta t}{2} \nabla_x V(x^{\ast \ast}, y^{\ast \ast}), \quad y^{n+1} = y^{\ast \ast}, \quad v^{n+1} = v^{\ast \ast} - \frac{\Delta t}{2} \nabla_y V(x^{\ast \ast}, y^{\ast \ast}).
\end{align*}

We provide two cases of proofs for our main theorem. One is in the simplified case when the initial wavefunction is a complex Gaussian and in this proof we only use the Wigner function; the second proof is for the general case, where the Wigner function loses its non-negativity, we are forced to turn to the Husimi functions. Both proofs employ the observable estimate of the QCMD \cref{thm:cts_egorov} and its consequence \cref{cor:egorov_wigner_husimi}, which bridge the quantum observables with a classical Hamiltonian flow, and the numerical errors of the St\"ormer–Verlet integrator of this classical counterpart.

\begin{thm}[Observable Convergence of the Strang Splitting]\label{thm:observable}
Let the observable $A = \op(a)$ be the Weyl quantization of a Schwartz function $a: \mathbb{R}^{2d} \to \mathbb{R}$, and the potential function $V$ be smooth with its derivatives of order $\geq$ 2 all bounded. Then the global error in the expectation value of $A$ can be estimated as
\begin{equation} \label{eqn:main_thm_ob}
|\braket{\psi^h_n |A| \psi^h_n} - \braket{\psi^h(T)|A| \psi^h(T)}| \leq C_T(\Delta t^2 + h^2),
\end{equation}
for some $C_T$ depends on $t_n = T$ and the derivatives of $a$ and $V$.
\end{thm}
\begin{proof}[Proof of the case with an initial complex gaussian]
From $t_n$ to $t_{n+1}$, because the construction of the scheme 
\begin{align*}
    \psi^{h}_{n+1} &= U_V(t_{n}+\Delta t/2, t_{n+1}) U_T(t_n, t_{n+1}) U_V (t_n, t_{n}+\Delta t/2)\psi^{h}_{n}
    \\
     &=  U_V(t_{n}+\Delta t/2, t_{n+1}) U_T(t_n, t_{n+1}) \psi^{h}_{\ast}
      = U_V(t_{n}+\Delta t/2, t_{n+1}) \psi^{h}_{\ast \ast}.
\end{align*}
and \cref{thm:cts_egorov}, one has
\begin{align*}
 \braket{\psi^h_{n+1} |A| \psi^h_{n+1}} 
 &= \braket{A}_{U_V(t_{n}+\Delta t/2, t_{n+1}) \psi^{h}_{\ast \ast}} \\
 &= \braket{\op(a \circ \Phi_V^{\Delta t/2}) }_{\psi^{h}_{\ast \ast}} + \Or (\Delta t h^2)
 \\
 &= \braket{\op(a \circ \Phi_V^{\Delta t/2} \circ  \Phi_T^{\Delta t} }_{\psi^{h}_{\ast}} + \Or (\Delta t h^2)
 \\
 &= \braket{\op(a \circ \Phi_V^{\Delta t/2} \circ  \Phi_T^{\Delta t} \circ \Phi_V^{\Delta t/2}) }_{\psi^{h}_{n}} + \Or (\Delta t h^2)
 \\
 & = \braket{\op(a \circ  \Phi^{\Delta t}_{SV}) }_{\psi^{h}_{n}} + \Or (\Delta t h^2)
,
\end{align*}
One can repeat this procedure for each time interval and obtain
\begin{align}
\braket{\psi^h_{n} |A| \psi^h_{n}} &= \braket{\op(a \circ  \left(\Phi^{\Delta t}_{SV} \right)^n ) }_{\psi^{h}_0} + \Or (T h^2)
\nonumber\\
&= \int_{\mathbb{R}^{2d}} \left(a\circ  \left(\Phi^{\Delta t}_{SV} \right)^n  \right)(q, p) w^h[{\psi^{h}_0}](q, p) \,dq\,dp +  \Or ( h^2),
\label{eqn:pf_wigner_scheme}
\end{align}
where in the last line we used the relationship of the Wigner function and Weyl quantization. Note that since $\psi_0^h$ is a complex Gaussian, its Wigner function is non-negative, i.e. $ w^h[{\psi_0^h}]\geq 0$ and 
\[
\int_{\mathbb{R}^{2d}}  w^h[{\psi^{h}_{0}}] \, dq \, dp = \norm{\psi^{h}_0}^2 = 1.
\]
On the other hand, applying \cref{cor:egorov_wigner_husimi} to the exact solution $\psi^h(t_n)$, we have
\begin{align}
\braket{\psi^h(t_n) |A| \psi^h(t_n)} 
= \int_{\mathbb{R}^{2d}} \left(a\circ \Phi^{t_n}  \right)(q, p)  w^h[{\psi^{h}_0}](q, p) \,dq\,dp +  \Or ( h^2),
\label{eqn:pf_wigner_cts}
\end{align}
and hence it suffices to compare the difference in the $a$ following the two different classical flows. The error bounds of the St\"omer-Verlet integrator of the classical Hamiltonian system read
\[
\norm{a\circ  \Phi^{t_n}  - a\circ (\Phi^{\Delta t}_{SV})^n }_{L^\infty(\mathbb{R}^{2d})} \leq C \Delta t^2, 
\]
for some constant $C$ depend on $t_n = T$. As an immediate consequence, one has
\begin{align*}
 \int_{\mathbb{R}^{2d}} \left(a\circ  \left(\Phi^{\Delta t}_{SV} \right)^n - a\circ \Phi^{t_n} \right)(q, p)  w^h[{\psi^{h}_0}](q, p) \,dq\,dp  \leq C \Delta t^2.
\end{align*}
Combining \eqref{eqn:pf_wigner_scheme} and \cref{eqn:pf_wigner_cts}, we arrived at the desired result.
\end{proof}

For general initial data, one no longer has the non-negativity of the Wigner function $  w^h[{\psi_0^h}] \geq 0$, but fortunately, the Husimi function is still non-negative. Therefore, we perform the proof using the Husimi function instead.
\begin{proof}[Proof of the general case]
Similar as the proof using the Wigner function, we have
\begin{align*}
\braket{\psi^h_{n} |A| \psi^h_{n}} &= \braket{\op(a \circ  \left(\Phi^{\Delta t}_{SV} \right)^n ) }_{\psi^{h}_0} + \Or (T h^2)
\\
&= \int_{\mathbb{R}^{2d}} \left(a\circ  \left(\Phi^{\Delta t}_{SV} \right)^n - \frac{h}{4}\Delta a\circ  \left(\Phi^{\Delta t}_{SV} \right)^n \right)(q, p) \sigma^h[{\psi^{h}_{0}}](q, p) \,dq\,dp +  \Or ( h^2),
\end{align*}
where the last line we used \cref{lmm:exp_husimi} to represent the expectation in the Husimi function. On the other hand, thanks to \cref{cor:egorov_wigner_husimi}, we also have a similar result of the exact solution
\begin{align*}
\braket{\psi^h(t_n) |A| \psi^h(t_n)} 
&= \int_{\mathbb{R}^{2d}} \left(a\circ  \Phi^{t_n}  - \frac{h}{4}\Delta a\circ \Phi^{t_n} \right)(q, p) 
\sigma^h[{\psi^{h}_0}](q, p) \,dq\,dp +  \Or ( h^2).
\end{align*}
Notice that $\sigma^h[{\psi_0}]$ is a non-negative function and satisfies
\[
\int_{\mathbb{R}^{2d}} \sigma^h[{\psi^{h}_0}] \, dq \, dp = \norm{\psi^{h}_{0}}^2 = 1,
\]
and hence besides the comparison of $a$ following the flows $\Phi^{t_n}$ and $(\Phi^{\Delta t}_{SV})^n$, we also need to compare that of $\Delta a$. Fortunately, the error bounds of the St\"omer-Verlet integrator \cite{HairerEtAl03} of the classical Hamiltonian system satisfy
\[
\lvert \Phi^{t_n}(q,p)  - (\Phi^{\Delta t}_{SV})^n(q,p) \rvert \leq C \Delta t^2,
\quad
\lvert D\Phi^{t_n}(q,p)  - D(\Phi^{\Delta t}_{SV})^n(q,p) \rvert \leq C \Delta t^2,
\]
for some constant $C$ depend on $t_n = T$ and independent for $n$ and $\Delta_t$ thanks to the symplecticity of the integrator. Thus,
\[
\norm{a\circ  \Phi^{t_n}  - a\circ (\Phi^{\Delta t}_{SV})^n }_{L^\infty(\mathbb{R}^{2d})} \leq C \Delta t^2, 
\]
and
\[
\norm{\Delta a\circ  \Phi^{t_n}  - \Delta a\circ (\Phi^{\Delta t}_{SV})^n }_{L^\infty(\mathbb{R}^{2d})} 
\leq \norm{D\left(\Delta a\right)}_{L^\infty(\mathbb{R}^{2d})} \norm{D\Phi^{t_n}  - D(\Phi^{\Delta t}_{SV})^n}_{L^\infty(\mathbb{R}^{2d})} 
\leq C \Delta t^2,
\]
Therefore, 
\[
\lvert \braket{\psi^h(t_n) |A| \psi^h(t_n)}  - \braket{\psi^h_{n} |A| \psi^h_{n}} \rvert \leq C_T (\Delta t^2 +h^2),
\]
as desired.
\end{proof}

Note that \cref{thm:observable} is not a uniform bound due to the presence of $h$ on the right hand side, but it drastically improves the result using a direct estimate of the wavefunction $\Delta t ^2 /h^2$ or $\Delta t^2 / h$ depending on the regularity of the initial condition (see, e.g., \cite{Bao:2002fy, Thalhammer1, Thalhammer2}). Nevertheless, we remark that it is possible to achieve an estimate uniform in $h$ combining the two kind of error bounds, namely, the error can be chosen as the minimum value between them
\begin{equation} \label{eqn:min_two_errors}
\min_{0<h\leq 1} \left\{ \Delta t^2 + h^2, \frac{\Delta t ^2}{h} \right\}
\end{equation}
and in the worst case scenario when $h = \Or(\Delta t^{2/3})$, one has the error bound of $\Or(\Delta t^{4/3})$. We summarize the result in the following corollary. 
\begin{cor}[A Uniform Observable Estimate of the Strang Splitting] \label{cor:uniform_ob}
Let the observable $A = \op(a)$ be the Weyl quantization of a Schwartz function $a: \mathbb{R}^{2d} \to \mathbb{R}$, and the initial wavefunction $\psi^h_0(x)$ is of the WKB type with bounded amplitude and phase. Under the assumptions of \cref{thm:observable},
we have the uniform-in-$h$ estimate of the global error in the expectation value of $A$ as
$$|\braket{\psi^h_n |A| \psi^h_n} - \braket{\psi^h(T)|A| \psi^h(T)}| \leq C_T \Delta t^{4/3},$$
for some $C_T$ depends on $t_n = T$, the functions $a$ and $V$ and the initial wavepacket $\psi^h_0$.
\end{cor}

\begin{rem}
\begin{enumerate}
    \item We remark that the additional assumption on the initial wavefunction comes from the condition of the result for the wavefunction error estimate. To be specific, for the WKB type initial condition $\psi^h_0$, one has $h^j\norm{\psi_0}_{H^j} \leq M_j$ with some constant $M_j$ for $j = 1, 2$, which is the assumption of the initial condition for the vector norm bound of such semiclassical problem~\cite{Thalhammer1,Thalhammer2}.
    
    \item More importantly, we remark that the convergence rate of $4/3$ is not sharp. The optimal rate is unknown in the literature even for the linear case. We conjecture that the optimal uniform-in-$h$ rate is $\Delta t^2$ at least for the linear case, and the improvement of the estimate is left for future work.
    
    \item We also clarify that these results are for the time-splitting (Trotter) errors with spatial degrees of freedom kept continuous. The spatial discretization error is not characterized here. In the numerical examples, we implement the spatial discretization via the Fourier-based pseudospectral discretization with $\Delta x = \Or(h)$, which has high accuracy (spectral accuracy). To prove a uniform error bound with spatial discretization is a more delicate issue, and there is so far no rigorous result in the literature even for the linear case, which remains as an important open question in this field. The description for the nonlinear case is beyond the scope of this paper. Nevertheless, it is worth mentioning that we have recently discovered an approach directly analyzing the observable error bounds in the spatially discrete setting via discrete microlocal analysis~\cite{Borns-Weil2022,ChristiansenZworski2010,DyatlovJezequel2021}, which is a work in progress by one of the authors~\cite{BornsweilFang2022}.
    
\end{enumerate}
\end{rem}

\section{Numerical Results} \label{sec:num_results}

In this section, we shall report on a few numerical examples, which illustrate our theoretical results, in particular the observable error bounds presented in \cref{thm:observable} and the errors of the wavefunction depicted in \eqref{eqn:min_two_errors}.

To this end, we choose an interaction potential of the form $V(x,y) = \sin(x^2 + y^2)$, and the initial condition as
\begin{equation}
    \psi^h_{0}=\mathcal{Z} \exp\left(-12.5(x+1)^{2} + \I 50 (x+1) \right),
\end{equation} 
where $\mathcal{Z}$ is the normalization factor such that $\norm{\psi^h_0}_{L^2} =1$. The computational domain is $[-\pi, \pi]$ with $\Delta x = 2 \pi h/32$. The numerical solutions are computed by the time-splitting methods till $T = 0.5$ and the trapezoidal rules are applied spatially for the quadrature that are known to be spectral accurate for periodic functions. We present the absolute error of the wavefunctions $\psi^h$ at $t_n = T$ are computed by 
\begin{equation*}
(\Delta x)^{1/2} \norm{\psi^h_n - \psi^h(t_n)}_2  \approx \norm{\psi^h_n(\cdot) - \psi^h(t_n, \cdot)}_{L^2},
\end{equation*}
and for the observable expectations are computed using absolute values between the reference and numerical values as defined on the left-hand-side of \eqref{eqn:main_thm_ob}. We use the numerical results computed by a very small step size $\Delta t = 10^{-5}$ as the reference solution. 

Besides the wavefunction, we consider the expectation $\langle A \rangle_{\psi^h(t, \cdot)}$ of the following observables:
\begin{itemize}
    \item The position operator $\hat x$, defined as multiplication by $x$, namely $\hat x \psi = x \psi$;
    \item The momentum operator $\hat p = - \I h \nabla_x$;
    \item The Gaussian observable (labelled as ``Gaussian" in the figures) defined as the multiplication by $e^{-4x^{2}}$;
    \item The Schwartz observable (labelled as ``xGaussian" in the figures) defined as the multiplication by $x e^{-4x^{2}}$;
    \item kinetic operator $\frac{1}{2} \hat p^2 = - \frac{1}{2} h^2 \Delta_x$.
\end{itemize}

\begin{figure}[h!]
        \centering
        \subfloat[Wavefunction $\psi^h$, $y$ and $v$ ]{
        \includegraphics[width=.48\textwidth]{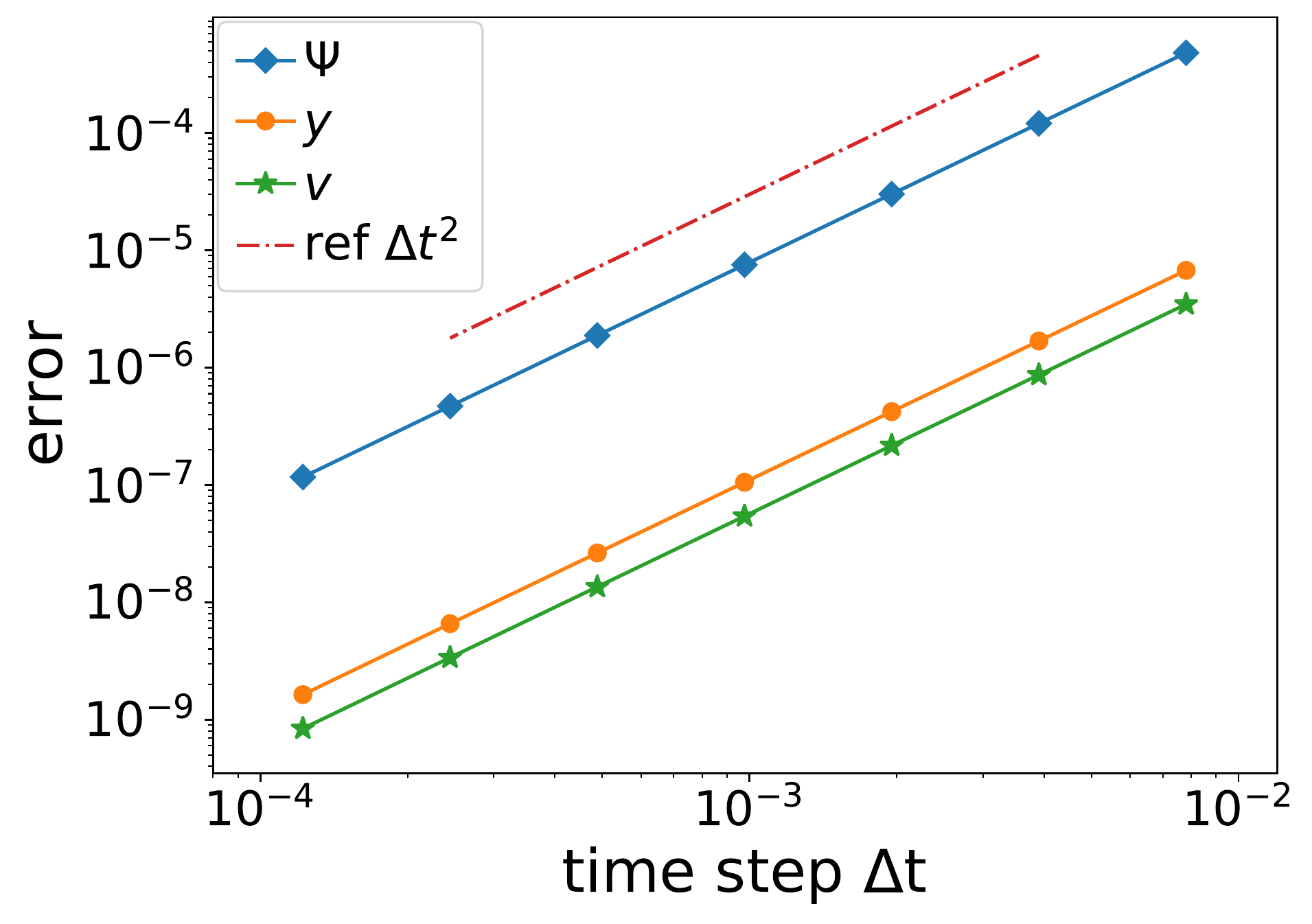}
        \label{fig:1a}}
        \subfloat[Observable expectations]{
        \includegraphics[width=.48\textwidth]{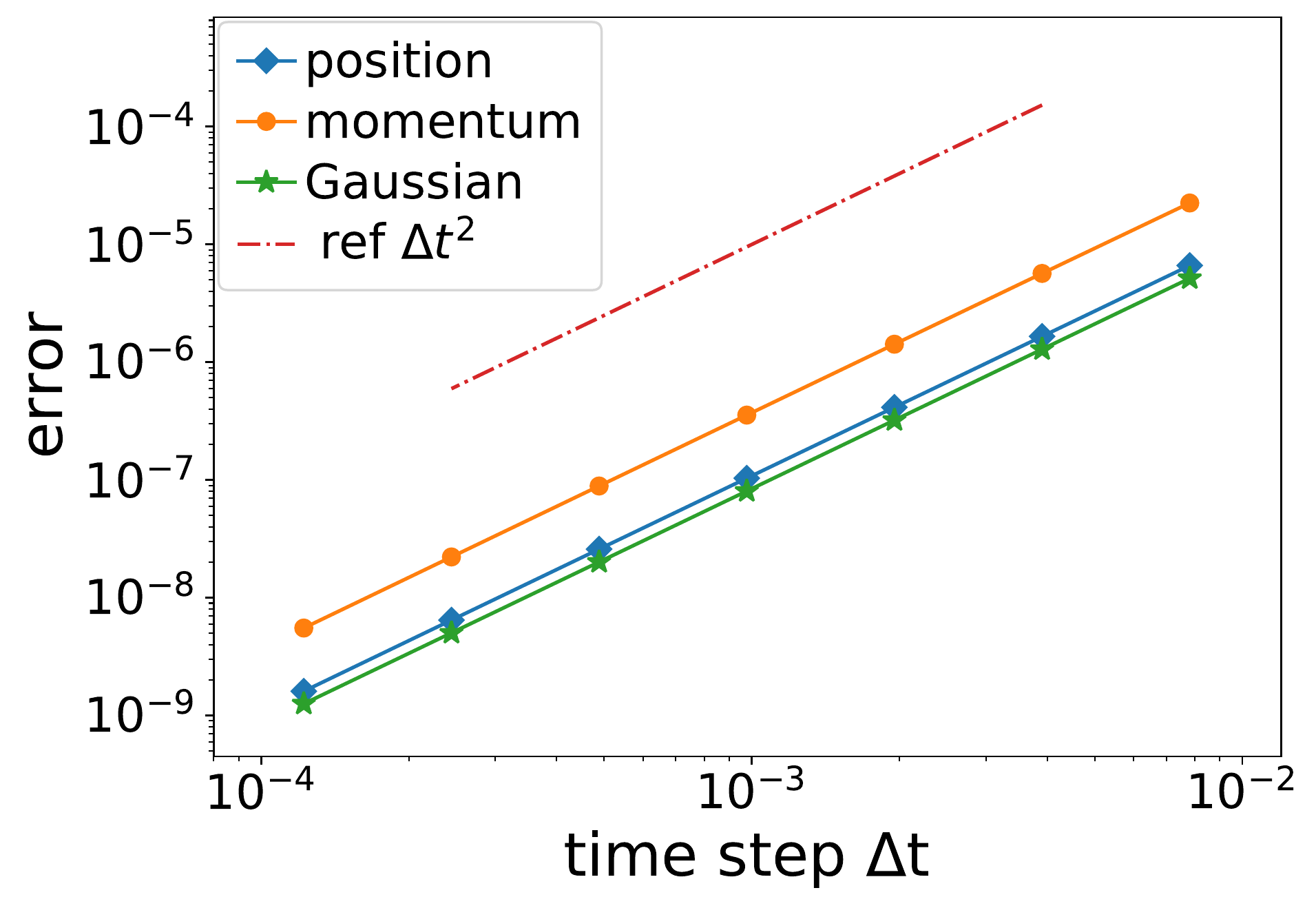}
        \label{fig:1b}}
        \caption{Errors versus various time steps. Left: Log-log plot of the errors of the wavefunctions $\psi$, $y$ and $v$. Right: Log-log plot of the errors of the position, momentum and Gaussian observable expectations. All errors are of second-order convergence in time.  }
        \label{fig:var_t_wav_ops}
    \end{figure}

To demonstrate \cref{thm:observable}, we present the numerical results of varying both $\Delta t$ and $h$. In the first test, we fix $h = 0.04$ and take the step sizes $\Delta t = 2^{-6}, 2^{-7}, \cdots, 2^{-11}$. The numerical errors of the wavefunction $\psi^h$ and the macroscopic quantities $y$ and $v$ are plotted in \cref{fig:1a}, while \cref{fig:1b} plots the errors in the expectations of the position, momentum and Gaussian observables. In can be seen that all errors are of second-order convergence in time, which confirms the scaling in $\Delta t$ as shown in \cref{eqn:main_thm_ob}. However, the errors in the wavefunction are of magnitude larger than those of the observable expectations and the macroscopic quantities $y$ and $v$. This implies that to capture the correct physical observables, one may choose larger step size comparing to the wavefunction.

\begin{figure}[h!]
        \centering
        \subfloat[Schwartz observables and Kinetic operators]{
        \includegraphics[width=.48\textwidth]{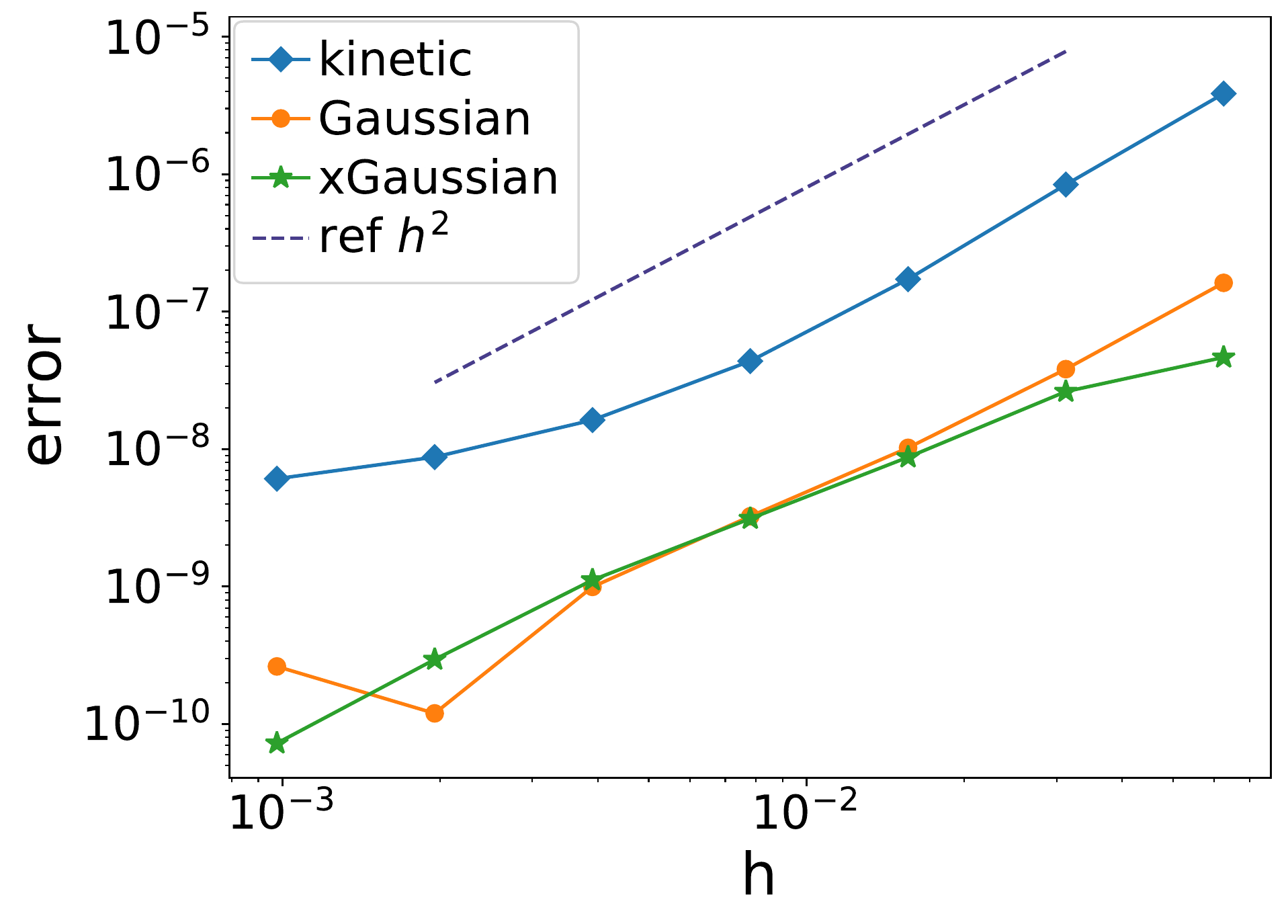}
        \label{fig:2a}}
        \subfloat[Position and momentum operators]{
        \includegraphics[width=.48\textwidth]{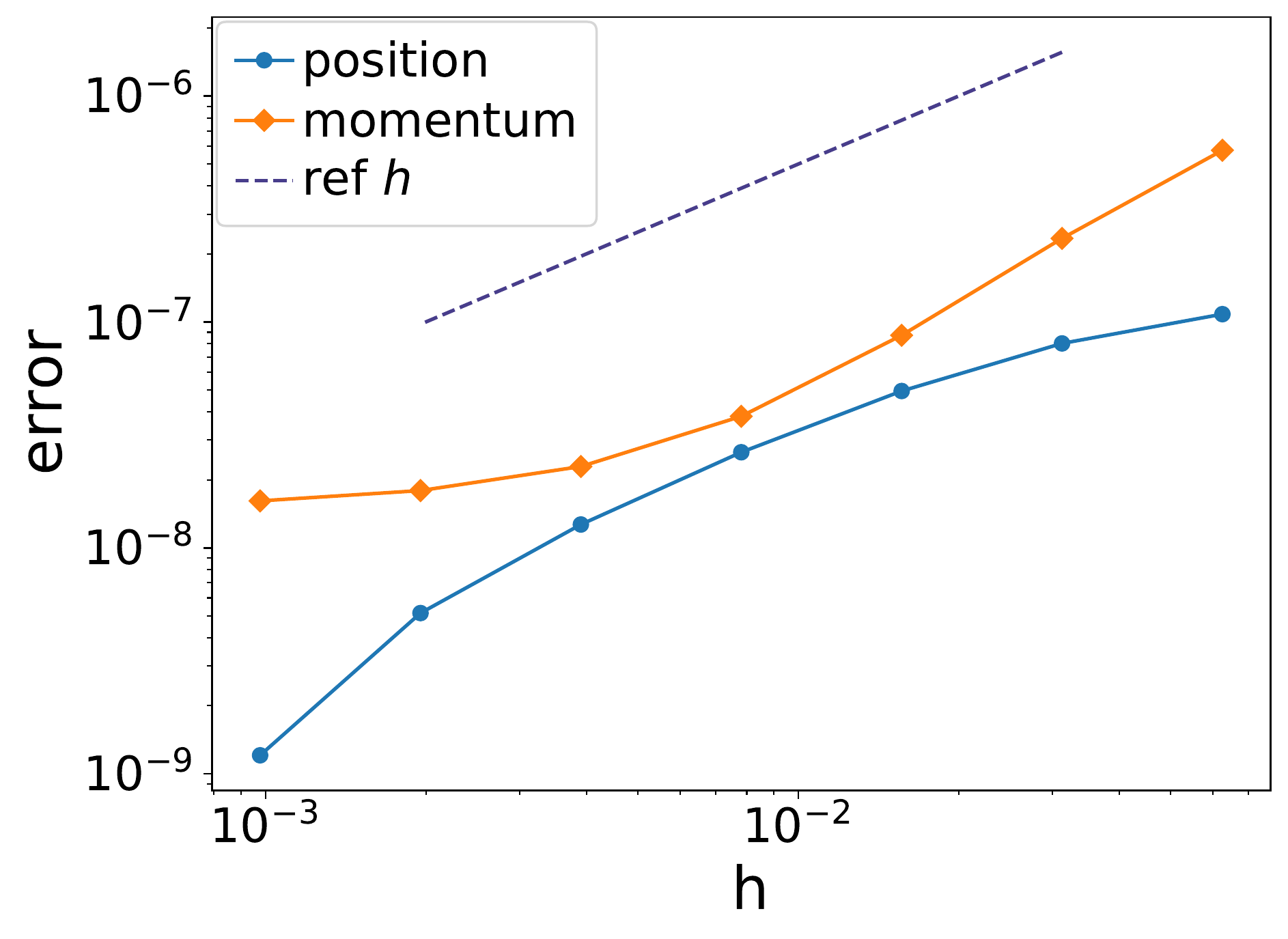}
        \label{fig:2b}} 
        \caption{Errors for various $h$. Left: Log-log plot of the errors of the expectations for the Schwartz observables and kinetic operator. Right: Log-log plot of the errors of the position and momentum operators. For Schwartz observables, the error scales as $h^2$, which matches our theoretical bounds as in \cref{eqn:main_thm_ob}. For non-Schwartz observables, a scaling of $h^\alpha$ (with $\alpha>0$) is observed.  }
        \label{fig:var_h_ops}
    \end{figure}

We then vary the size of the semiclassical parameter $h$. In this test, we fix $\Delta t = 0.001$ and take $h = 2^{-4}, 2^{-5}, \cdots, 2^{-10}$. The errors of the observable expectations are recorded in \cref{fig:var_h_ops}. It can be seen in \cref{fig:2a} that the errors of the Gaussian and Schwartz observables are of second order in $h$, which agrees with \cref{eqn:main_thm_ob}. Besides the Schwartz observables, we also test other non-Schwartz observables such as the position, momentum and kinetic operators. Note that although in these  cases, our main theorem can no longer be applied, we still observe additive error bounds in $h$. In particular, the errors of the expectations for the position and momentum operator scale as $\Or(h)$ as shown in \cref{fig:2b} and that for the position operator scales as $\Or(h^2)$. It is worth pointing out that we define the Weyl quantization \eqref{eqn:def_weyl} for Schwartz functions $a$ and the lemma on expectation \cref{lmm:exp_husimi} also requires $a$ to be Schwartz, and hence a proof of rigorous observable error bounds for such non-Schwartz observables or polynomials of $\hat x$ and $\hat p$ is beyond the current framework, which is left as an interesting future direction.

\begin{figure}[h!]
        \centering
        \includegraphics[width=.48\textwidth]{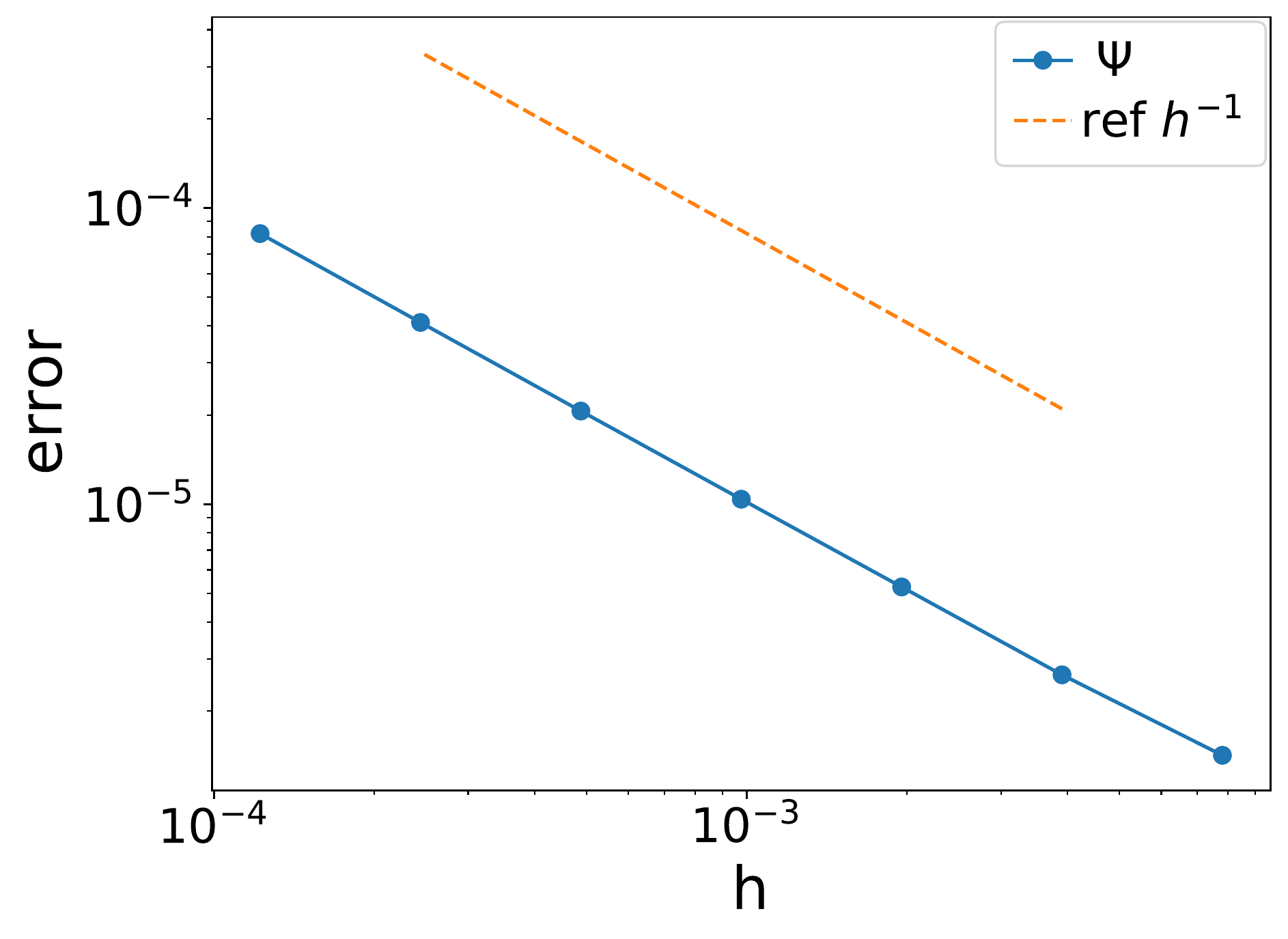}
        \caption{Log-log plot of the errors of the wavefunction $\psi^h$.}
        \label{fig:var_h_wav}
    \end{figure}

Finally, we present the error of the wavefunction when varying $h$. \cref{fig:var_h_wav} plots the errors of the wavefunction $\psi^h$ for $h = 2^{-7}, 2^{-8}, \cdots 2^{-13}$. It can be seen the scaling is $h^{-1}$ as given in \eqref{eqn:min_two_errors}, which is used to generate the uniform error bound. This $h^{-1}$ scaling is interesting in the following two perspectives. First, it helps to demonstrate the uniform error bounds. To be specific, note that the uniform-in-$h$ observable error bounds in \cref{cor:uniform_ob} is not tight and hence difficult to be directly numerically verified. Here we verify both terms in \eqref{eqn:min_two_errors} numerically instead. Furthermore, the error grows as $h$ increases which is the opposite of the observable and this further illustrates that capturing the correct physical observables allows one to take an $\Or(1)$ time step, much larger than that of the wavefunction calculations.

\section{Conclusion and Discussion}

In this work, we investigate the distance between the QCMD and its semiclassical limit by a careful estimate of the observable expectations. This observable results on the continuous level are then rewritten using both the Wigner and Husimi functions. We then discuss the time-splitting strategies for this quantum-classical mixed dynamics, and prove the observable error bounds of the Strang splitting using the continuous-level result. In particular, the observable error bounds present an additive scaling as $\Delta t^2 + h^2$, which drastically improves the estimate comparing to a direct estimation via the errors in the wavefunction that scales as $\Delta t^2/h$. Finally, a uniform observable error bound is provided by combining both estimating strategies.

Possible future directions branch into three different streams. One is to go beyond the Strang splitting and consider higher order time-splitting strategies such as Suzuki construction \cite{Suzuki91}, Yoshida triple jump \cite{Yoshida90}, symmetric Zassenhaus splitting \cite{IserlesKropielnickaSingh18}, symmetric-conjugate composition \cite{blanes2021symmetricconjugate}, for this nonlinear system. The other is to investigate such observable estimates for other linear and nonlinear systems with their corresponding numerical schemes, such as the surface hopping methods with matrix potentials \cite{ChaiJinLi13,ChaiJinLi15,LuZhou18,FangLu2018}, the time-dependent self-consistent field equations \cite{Jin:2017bh, tdscf1, tdscf3, tdscf2, tdscf4}, and other mean-field systems. We point out that the QCMD preserves its total energy, and hence in the semiclassical limit an autonomous Hamiltonian flow is arrived at. However, when the total energy is not preserved or the symbol $\mathcal{H}_0$ of the total energy is time dependent, one gets a non-autonomous ODE system instead. Even though the Egorov's theorem \cref{thm:cts_egorov} on the continuous level could still hold, the analysis for the time-splitting schemes may require a correspondence between the time-splitting strategy with some symplectic integrator for non-autonomous Hamiltonian system, which could bring some challenges to the analysis. Nevertheless, such non-autonomous cases would be an interesting future direction. The last is to consider non-Schwartz functions, such as polynomials of the position and momentum operators, which is left for future study.

Our paper focuses on the observable error bounds in terms of the time discretization. On a separate matter, to estabilish a uniform error bound of the observables in the spatially discrete setting is an important open problem even for the linear cases. Though beyond the scope of this paper, we have recently discovered a promising approach directly analyzing the observable error bounds with spatial discretizations via discrete microlocal analysis~\cite{Borns-Weil2022,ChristiansenZworski2010,DyatlovJezequel2021,BornsweilFang2022}, which is a work in progress by one of the authors. As for implimentations of the high-dimensional PDE grids, recent advances of quantum algorithms for unbounded Hamiltonian simulation~\cite{BabbushBerryEtAl2019,AnEtAl2021,SuBerryEtAl2021,TongAlbertEtAl2021,AnFangLin2022,ChildsLengEtAl2022} provide new prospect for efficent implementations and the combination of the observable error bounds and quantum implementations is an interesting future direction.

\section*{Acknowledgment}

DF is supported by NSF Quantum Leap Challenge Institute (QLCI) program under Grant No. OMA-2016245, NSF DMS-2208416, and a grant from the Simons Foundation under Award No. 825053. The authors also acknowledge the hospitality of Kavli Institute for Theoretical Physics (KITP) supported by the National Science Foundation under Grant No. NSF PHY-1748958.

\bibliographystyle{abbrv}
\bibliography{qcmd}

\end{document}